\documentclass[leqno]{amsart}
\usepackage{amsthm}
\usepackage{amssymb}
\usepackage{esint}
\usepackage[unicode=true,pdfusetitle,
 bookmarks=true,bookmarksnumbered=false,bookmarksopen=false,
 breaklinks=false,pdfborder={0 0 1},backref=false,colorlinks=false]
 {hyperref}

\makeatletter
\allowdisplaybreaks

\usepackage[numbers,sort&compress,square]{natbib}

\makeatother

\begin{document}

\title{Identities of some special mixed-type polynomials }

\author{Dae San Kim}

\address{\noindent Department of Mathematics, Sogang University, Seoul 121-742,
Republic of Korea}

\email{dskim@sogang.ac.kr}

\author{Taekyun Kim}

\address{\noindent Department of Mathematics, Kwangwoon University, Seoul
139-701, Republic of Korea}

\email{tkkim@kw.ac.kr}


\begin{abstract}
In this paper, we consider various speical mixed-type polynomials which are
related to Bernoulli, Euler, Changhee and Daehee polynomials. From
those polynomials, we derive some interesting and new identities.
\end{abstract}
\maketitle

\global\long\def\Zp{\mathbb{Z}_{p}}

\section{Introduction}

Let $p$ be a fixed odd prime number. Throughout this paper, $\Zp$,
$\mathbb{Q}_{p}$ and $\mathbb{C}_{p}$ will denote the ring of $p$-adic
integers, the field of $p$-adic rational numbers and the completion
of the algebraic closure of $\mathbb{Q}_{p}$. Let $\nu_{p}$ be the
normalized exponential valuation of $\mathbb{C}_{p}$ with $\left|p\right|_{p}=\frac{1}{p}=p^{-\nu_{p}\left(p\right)}$.
Let $UD\left(\Zp\right)$ be the space of uniformly differentiable
functions on $\Zp$.

For $f\in UD\left(\Zp\right)$, the bosonic $p$-adic integral is
given by
\begin{equation}
I_{0}\left(f\right)=\int_{\Zp}f\left(x\right)d\mu_{0}\left(x\right)=\lim_{N\rightarrow\infty}\frac{1}{p^{N}}\sum_{x=0}^{p^{N}-1}f\left(x\right),\quad\left(\mbox{see \cite{key-11}}\right),\label{eq:1}
\end{equation}
and the fermionic $p$-adic integral on $\Zp$ is defined by Kim to
be
\[
I_{-1}\left(f\right)=\int_{\Zp}f\left(x\right)d\mu_{-1}\left(x\right)=\lim_{N\rightarrow\infty}\sum_{x=0}^{p^{N}-1}f\left(x\right)\left(-1\right)^{x},\quad\left(\mbox{see \cite{key-10}}\right).
\]

In \cite{key-6,key-7}, the higher-order Daehee polynomials are defined
by
\begin{equation}
\left(\frac{\log\left(1+t\right)}{t}\right)^{r}\left(1+t\right)^{x}=\sum_{n=0}^{\infty}D_{n}^{\left(r\right)}\left(x\right)\frac{t^{n}}{n!},\quad\left(r\in\mathbb{N}\right).\label{eq:2}
\end{equation}

When $x=0$, $D_{n}^{\left(r\right)}=D_{n}^{\left(r\right)}\left(0\right)$
are called the Daehee numbers of order $r$.

When $r=1$, $D_{n}^{\left(1\right)}\left(x\right)=D_{n}\left(x\right)$
are called the Daehee polynomials (see \cite{key-6}).

As is known, the Changhee polynomials of order $s$$\left(\in\mathbb{N}\right)$
are defined by the generating function to be
\begin{equation}
\left(\frac{2}{t+2}\right)^{s}\left(1+t\right)^{x}=\sum_{n=0}^{\infty}Ch_{n}^{\left(s\right)}\left(x\right)\frac{t^{n}}{n!},\quad\left(\mbox{see \cite{key-8}}\right).\label{eq:3}
\end{equation}

When $x=0$, $Ch_{n}^{\left(s\right)}=Ch_{n}^{\left(s\right)}\left(0\right)$
are called the Changhee numbers of order $s$.

For $s=1$, $Ch_{n}^{\left(1\right)}\left(x\right)=Ch_{n}\left(x\right)$
are called the Changhee polynomials.

The Bernoulli polynomials of order $r\in\mathbb{N}$ are given by

\begin{equation}
\left(\frac{t}{e^{t}-1}\right)^{r}e^{xt}=\sum_{n=0}^{\infty}B_{n}^{\left(r\right)}\left(x\right)\frac{t^{n}}{n!},\quad\left(\mbox{see \cite{key-11,key-12,key-13,key-14,key-15,key-16,key-17,key-18,key-19,key-20,key-21}}\right).\label{eq:4}
\end{equation}

When $x=0$, $B_{n}^{\left(r\right)}=B_{n}^{\left(r\right)}\left(0\right)$
are called the Bernoulli numbers of order $r$.

For $r=1$, $B_{n}^{\left(1\right)}\left(x\right)=B_{n}\left(x\right)$
are called the ordinary Bernoulli polynomials.

We recall that the Euler polynomials of order $r$ are defined by
the generating function to be

\begin{equation}
\left(\frac{2}{e^{t}+1}\right)^{r}e^{xt}=\sum_{n=0}^{\infty}E_{n}^{\left(r\right)}\left(x\right)\frac{t^{n}}{n!},\quad\left(\mbox{see \cite{key-1,key-2,key-3,key-4,key-5,key-6,key-7,key-8,key-9,key-10,key-11,key-12}}\right).\label{eq:5}
\end{equation}

When $x=0$, $E_{n}^{\left(r\right)}=E_{n}^{\left(r\right)}\left(0\right)$
are called the Euler numbers of order $r$.

For $r=1$, $E_{n}^{\left(1\right)}\left(x\right)=E_{n}\left(x\right)$
are called the ordinary Euler polynomials.

Finally, the Cauchy polynomials of the first kind of order $r$ are given by
\begin{equation}
\left(\frac{t}{\log\left(1+t\right)}\right)^{r}\left(1+t\right)^{x}=\sum_{n=0}^{\infty}C_{n}^{\left(r\right)}\left(x\right)\frac{t^{n}}{n!},\quad\left(\mbox{see \cite{key-3,key-9}}\right).\label{eq:6}
\end{equation}

When $x=0$, $C_{n}^{\left(r\right)}=C_{n}^{\left(r\right)}\left(0\right)$
are called the Cauchy numbers of the first kind of order $r$.

For $r=1$, $C_{n}^{\left(1\right)}\left(x\right)=C_{n}\left(x\right)$
are called the ordinary Cauchy polynomials of the first kind (see \cite{key-3}).

From (\ref{eq:1}) and (\ref{eq:2}), we have
\begin{equation}
I_{0}\left(f_{1}\right)-I_{0}\left(f\right)=f^{\prime}\left(0\right)\label{eq:7}
\end{equation}

and
\begin{equation}
I_{-1}\left(f_{1}\right)=-I_{-1}\left(f\right)+2f\left(0\right),\label{eq:8}
\end{equation}
where $f_{1}\left(x\right)=f\left(x+1\right)$ (see \cite{key-10,key-12}).

In this paper, we consider several special polynomials which are
derived from the bosonoic or fermionic $p$-adic integral on $\Zp$.

Finally, we give some relation or identities of those polynomials.

\section{Some special mixed-type polynomials }

In this section, we assume that $t\in\mathbb{C}_{p}$ with $\left|t\right|_{p}<p^{-\frac{1}{p-1}}$.
From (\ref{eq:7}), we can derive the following equation :
\begin{align}
 & \int_{\Zp}\cdots\int_{\Zp}\left(1+t\right)^{\left(x_{1}+\cdots+x_{r}+x\right)}d\mu_{0}\left(x_{1}\right)\cdots d\mu_{0}\left(x_{r}\right)\label{eq:9}\\
= & \left(\frac{\log\left(1+t\right)}{t}\right)^{r}\left(1+t\right)^{x}=\left(\frac{\log\left(1+t\right)}{e^{\log\left(1+t\right)}-1}\right)^{r}e^{x\log(1+t)}\nonumber \\
=& \sum_{m=0}^{\infty}B_{m}^{(r)}(x)\frac{(\log(1+t))^{m}}{m!}=\sum_{m=0}^{\infty}B_{m}^{(r)}(x)\sum_{n=m}^{\infty}S_1(n,m)\frac{t^{n}}{n!}\nonumber\\
= & \sum_{n=0}^{\infty}\left(\sum_{m=0}^{n}B_{m}^{\left(r\right)}\left(x\right)S_{1}\left(n,m\right)\right)\frac{t^{m}}{n!},\nonumber
\end{align}
and
\begin{equation}
\left(\frac{\log\left(1+t\right)}{t}\right)^{r}\left(1+t\right)^{x}=\sum_{n=0}^{\infty}D_{n}^{\left(r\right)}\left(x\right)\frac{t^{n}}{n!}.\label{eq:10}
\end{equation}

Therefore, by (\ref{eq:9}) and (\ref{eq:10}), we obtain the following
equation :
\begin{align}
 & \int_{\Zp}\cdots\int_{\Zp}\dbinom{x_{1}+\cdots+x_{r}+x}{n}d\mu_{0}\left(x_{1}\right)\cdots d\mu_{0}\left(x_{r}\right)\label{eq:11}\\
= & \frac{D_{n}^{\left(r\right)}\left(x\right)}{n!}=\frac{1}{n!}\sum_{m=0}^{n}B_{m}^{\left(r\right)}\left(x\right)S_{1}\left(n,m\right)\nonumber
\end{align}
where $S_{1}\left(n,m\right)$ is the Stirling number of the first
kind.

From (\ref{eq:8}), we have
\begin{align}
 & \int_{\Zp}\cdots\int_{\Zp}\left(1+t\right)^{x_{1}+\cdots+x_{r}+x}d\mu_{-1}\left(x_{1}\right)\cdots d\mu_{-1}\left(x_{r}\right)\label{eq:12}\\
= & \left(\frac{2}{t+2}\right)^{r}\left(1+t\right)^{x}=\left(\frac{2}{e^{\log\left(1+t\right)}+1}\right)^{r}e^{x\log\left(1+t\right)}\nonumber \\
= & \sum_{m=0}^{\infty}E_{m}^{\left(r\right)}\left(x\right)\frac{\left(\log\left(1+t\right)\right)^{m}}{m!}=\sum_{m=0}^{\infty}E_{m}^{\left(r\right)}\left(x\right)\sum_{n=m}^{\infty}S_{1}\left(n,m\right)\frac{t^{n}}{n!}\nonumber \\
= & \sum_{n=0}^{\infty}\left\{ \sum_{m=0}^{n}E_{m}^{\left(r\right)}\left(x\right)S_{1}\left(n,m\right)\right\} \frac{t^{n}}{n!}\nonumber
\end{align}
and
\begin{equation}
\left(\frac{2}{t+2}\right)^{r}\left(1+t\right)^{x}=\sum_{n=0}^{\infty}Ch_{n}^{\left(r\right)}\left(x\right)\frac{t^{n}}{n!}.\label{eq:13}
\end{equation}

From (\ref{eq:12}) and (\ref{eq:13})
\begin{align}
 & \int_{\Zp}\cdots\int_{\Zp}\dbinom{x_{1}+\cdots+x_{r}+x}{n}d\mu_{-1}\left(x_{1}\right)\cdots d\mu_{-1}\left(x_{r}\right)\label{eq:14}\\
= & \frac{Ch_{n}^{\left(r\right)}\left(x\right)}{n!}=\frac{1}{n!}\sum_{m=0}^{n}E_{m}^{\left(r\right)}\left(x\right)S_{1}\left(n,m\right).\nonumber
\end{align}

Note that
\begin{align}
\left(1+t\right)^{x} & =\left(\frac{t}{\log\left(1+t\right)}\right)^{r}\left(1+t\right)^{x}\left(\frac{\log\left(1+t\right)}{t}\right)^{r}\label{eq:15}\\
 & =\left(\sum_{l=0}^{\infty}C_{l}^{\left(r\right)}\left(x\right)\frac{t^{l}}{l!}\right)\left(\sum_{m=0}^{\infty}D_{m}^{\left(r\right)}\frac{t^{m}}{m!}\right)\nonumber \\
 & =\sum_{n=0}^{\infty}\left(\sum_{l=0}^{n}\dbinom{n}{l}C_{l}^{\left(r\right)}\left(x\right)D_{n-l}^{\left(r\right)}\right)\frac{t^{n}}{n!}\nonumber
\end{align}
and
\begin{equation}
\left(1+t\right)^{x}=\sum_{n=0}^{\infty}\left(x\right)_{n}\frac{t^{n}}{n!}.\label{eq:16}
\end{equation}

From (\ref{eq:15}) and (\ref{eq:16}), we have
\begin{align}
\left(x\right)_{n} & =\sum_{l=0}^{n}\dbinom{n}{l}C_{l}^{\left(r\right)}\left(x\right)D_{n-l}^{\left(r\right)}\label{eq:17}\\
 & =\sum_{l=0}^{n}\dbinom{n}{l}D_{n-l}^{\left(r\right)}\left(x\right)C_{l}^{\left(r\right)}.\nonumber
\end{align}

That is,
\[
\dbinom{x}{n}=\frac{1}{n!}\sum_{l=0}^{n}\dbinom{n}{l}C_{l}^{\left(r\right)}\left(x\right)D_{n-l}^{\left(r\right)}.
\]

Let us consider the Bernoulli-Euler mixed-type polynomials of order $(r,s)$ as follows
:
\begin{equation}
BE_{n}^{\left(r,s\right)}\left(x\right)=\int_{\Zp}\cdots\int_{\Zp}E_{n}^{\left(s\right)}\left(x+y_{1}+\cdots+y_{r}\right)d\mu_{0}\left(y_{1}\right)\cdots d\mu_{0}\left(y_{r}\right).\label{eq:18}
\end{equation}

Then, we can find the generating function of $BE_{n}^{\left(r,s\right)}\left(x\right)$
as follows :
\begin{align}
 & \sum_{n=0}^{\infty}BE_{n}^{\left(r,s\right)}\left(x\right)\frac{t^{n}}{n!}\label{eq:19}\\
= & \int_{\Zp}\cdots\int_{\Zp}\sum_{n=0}^{\infty}E_{n}^{\left(s\right)}\left(x+y_{1}+\cdots+y_{r}\right)\frac{t^{n}}{n!}d\mu_{0}\left(y_{1}\right)\cdots d\mu_{0}\left(y_{r}\right)\nonumber \\
= & \left(\frac{2}{e^{t}+1}\right)^{s}\int_{\Zp}\cdots\int_{\Zp}e^{\left(x+y_{1}+\cdots+y_{r}\right)t}d\mu_{0}\left(y_{1}\right)\cdots d\mu_{0}\left(y_{r}\right)\nonumber \\
= & \left(\frac{2}{e^{t}+1}\right)^{s}\left(\frac{t}{e^{t}-1}\right)^{r}e^{xt}.\nonumber
\end{align}

Note that
\begin{align}
\left(\frac{2}{e^{t}+1}\right)^{s}\left(\frac{t}{e^{t}-1}\right)^{r}e^{xt} & =\left(\sum_{l=0}^{\infty}E_{l}^{\left(s\right)}\frac{t^{l}}{l!}\right)\left(\sum_{m=0}^{\infty}B_{m}^{\left(r\right)}\left(x\right)\frac{t^{m}}{m!}\right)\label{eq:20}\\
 & =\sum_{n=0}^{\infty}\left(\sum_{l=0}^{n}\dbinom{n}{l}E_{l}^{\left(s\right)}B_{n-l}^{\left(r\right)}\left(x\right)\right)\frac{t^{n}}{n!}.\nonumber
\end{align}

From (\ref{eq:19}) and (\ref{eq:20}), we have
\begin{equation}
BE_{n}^{\left(r,s\right)}\left(x\right)=\sum_{l=0}^{n}\dbinom{n}{l}E_{l}^{\left(s\right)}B_{n-l}^{\left(r\right)}\left(x\right).\label{eq:21}
\end{equation}

By replacing $t$ by $\log\left(1+t\right)$, we get
\begin{align}
 & \sum_{n=0}^{\infty}BE_{n}^{\left(r,s\right)}\left(x\right)\frac{\left(\log\left(1+t\right)\right)^{n}}{n!}\label{eq:22}\\
= & \left(\frac{2}{t+2}\right)^{s}\left(\frac{\log\left(1+t\right)}{t}\right)^{r}\left(1+t\right)^{x}\nonumber \\
= & \left(\sum_{l=0}^{\infty}Ch_{l}^{\left(s\right)}\frac{t^{l}}{l!}\right)\left(\sum_{m=0}^{\infty}D_{m}^{\left(r\right)}\left(x\right)\frac{t^{m}}{m!}\right)\nonumber \\
= & \sum_{n=0}^{\infty}\left\{ \sum_{m=0}^{n}\dbinom{n}{m}D_{m}^{\left(r\right)}\left(x\right)Ch_{n-m}^{\left(s\right)}\right\} \frac{t^{n}}{n!},\nonumber
\end{align}
and
\begin{align}
\sum_{m=0}^{\infty}BE_{m}^{\left(r,s\right)}\left(x\right)\frac{\left(\log\left(1+t\right)\right)^{m}}{m!} & =\sum_{m=0}^{\infty}BE_{m}^{\left(r,s\right)}\left(x\right)\sum_{n=m}^{\infty}S_{1}\left(n,m\right)\frac{t^{n}}{n!}\label{eq:23}\\
 & =\sum_{n=0}^{\infty}\left\{ \sum_{m=0}^{n}BE_{m}^{\left(r,s\right)}\left(x\right)S_{1}\left(n,m\right)\right\} \frac{t^{n}}{n!}.\nonumber
\end{align}

Therefore, by (\ref{eq:22}) and (\ref{eq:23}), we obtain the following
equation :
\begin{equation}
\sum_{m=0}^{n}\dbinom{n}{m}D_{m}^{\left(r\right)}\left(x\right)Ch_{n-m}^{\left(s\right)}=\sum_{m=0}^{n}BE_{m}^{\left(r,s\right)}\left(x\right)S_{1}\left(n,m\right).\label{eq:24}
\end{equation}

Let us consider the Daehee-Changhee mixed-type polynomials of order $(r,s)$ as follows
:

\begin{equation}
DC_{n}^{\left(r,s\right)}\left(x\right)=\int_{\Zp}\cdots\int_{\Zp}D_{n}^{\left(r\right)}\left(x+y_{1}+\cdots+y_{s}\right)d\mu_{-1}\left(y_{1}\right)\cdots d\mu_{-1}\left(y_{r}\right),\label{eq:25}
\end{equation}
where $n\ge0$.

From (\ref{eq:25}), we can derive the generating function of $DC_{n}^{\left(r,s\right)}\left(x\right)$
as follows :
\begin{align}
 & \sum_{n=0}^{\infty}DC_{n}^{\left(r,s\right)}\left(x\right)\frac{t^{n}}{n!}\label{eq:26}\\
= & \int_{\Zp}\cdots\int_{\Zp}\sum_{n=0}^{\infty}D_{n}^{\left(r\right)}\left(x+y_{1}+\cdots+y_{s}\right)\frac{t^{n}}{n!}d\mu_{-1}\left(y_{1}\right)\cdots d\mu_{-1}\left(y_{s}\right)\nonumber \\
= & \left(\frac{\log\left(1+t\right)}{t}\right)^{r}\int_{\Zp}\cdots\int_{\Zp}\left(1+t\right)^{x+y_{1}+\cdots+y_{s}}d\mu_{-1}\left(y_{1}\right)\cdots d\mu_{-1}\left(y_{s}\right)\nonumber \\
= & \left(\frac{\log\left(1+t\right)}{t}\right)^{r}\left(\frac{2}{t+2}\right)^{s}\left(1+t\right)^{x}.\nonumber
\end{align}

We observe that
\begin{align}
\left(\frac{2}{t+2}\right)^{s}\left(\frac{\log\left(1+t\right)}{t}\right)^{r}\left(1+t\right)^{x} & =\left(\sum_{l=0}^{\infty}Ch_{l}^{\left(s\right)}\frac{t^{l}}{l!}\right)\left(\sum_{m=0}^{\infty}D_{m}^{\left(r\right)}\left(x\right)\frac{t^{m}}{m!}\right)\label{eq:27}\\
 & =\sum_{n=0}^{\infty}\left\{ \sum_{m=0}^{n}\dbinom{n}{m}D_{m}^{\left(r\right)}\left(x\right)Ch_{n-m}^{\left(s\right)}\right\} \frac{t^{n}}{n!}.\nonumber
\end{align}

From (\ref{eq:26}) and (\ref{eq:27}), we have
\begin{equation}
DC_{n}^{\left(r,s\right)}\left(x\right)=\sum_{m=0}^{n}\dbinom{n}{m}D_{m}^{\left(r\right)}\left(x\right)Ch_{n-m}^{\left(r\right)},\label{eq:28}
\end{equation}
 where $n\ge0$, $r,s\in\mathbb{N}$.

Now, we define the Cauchy-Daehee mixed-type polynomials of order $(r,s)$ as follows
:
\begin{equation}
CD_{n}^{\left(r,s\right)}\left(x\right)=\int_{\Zp}\cdots\int_{\Zp}C_{n}^{\left(r\right)}\left(x+y_{1}+\cdots+y_{s}\right)d\mu_{0}\left(y_{1}\right)\cdots d\mu_{0}\left(y_{r}\right).\label{eq:29}
\end{equation}

From (\ref{eq:29}), we can derive the generating function of $CD_{n}^{\left(r,s\right)}\left(x\right)$
as follows :
\begin{align}
 & \sum_{n=0}^{\infty}CD_{n}^{\left(r,s\right)}\left(x\right)\frac{t^{n}}{n!}\label{eq:30}\\
= & \int_{\Zp}\cdots\int_{\Zp}\sum_{n=0}^{\infty}C_{n}^{\left(r\right)}\left(x+y_{1}+\cdots+y_{s}\right)\frac{t^{n}}{n!}d\mu_{0}\left(y_{1}\right)\cdots d\mu_{0}\left(y_{s}\right)\nonumber \\
= & \left(\frac{t}{\log\left(1+t\right)}\right)^{r}\int_{\Zp}\cdots\int_{\Zp}\left(1+t\right)^{x+y_{1}+\cdots+y_{s}}d\mu_{0}\left(y_{1}\right)\cdots d\mu_{0}\left(y_{s}\right)\nonumber \\
= & \left(\frac{t}{\log\left(1+t\right)}\right)^{r}\left(\frac{\log\left(1+t\right)}{t}\right)^{s}\left(1+t\right)^{x}.\nonumber \\
= & \begin{cases}
\sum_{n=0}^{\infty}C_{n}^{\left(r-s\right)}\left(x\right)\frac{t^{n}}{n!} & \mbox{if }r>s\\
\sum_{n=0}^{\infty}D_{n}^{\left(s-r\right)}\left(x\right)\frac{t^{n}}{n!} & \mbox{if }r<s\\
\sum_{n=0}^{\infty}\left(x\right)_{n}\frac{t^{n}}{n!} & \mbox{if }r=s.
\end{cases}\nonumber
\end{align}

Thus, by (\ref{eq:30}), we get
\begin{equation}
CD_{n}^{\left(r,s\right)}\left(x\right)=\begin{cases}
C_{n}^{\left(r-s\right)}\left(x\right) & \mbox{if }r>s\\
D_{n}^{\left(s-r\right)}\left(x\right) & \mbox{if }r<s\\
\left(x\right)_{n} & \mbox{if }r=s
\end{cases}\label{eq:31}
\end{equation}
where $n\ge0$.

By replacing $t$ by $e^{t}-1$ in (\ref{eq:26}), we get
\begin{align}
\sum_{n=0}^{\infty}DC_{n}^{\left(r,s\right)}\left(x\right)\frac{\left(e^{t}-1\right)^{n}}{n!} & =\left(\frac{t}{e^{t}-1}\right)^{r}e^{xt}\left(\frac{2}{e^{t}+1}\right)^{s}\label{eq:32}\\
 & =\left(\sum_{n=0}^{\infty}B_{l}^{\left(r\right)}\left(x\right)\frac{t^{l}}{l!}\right)\left(\sum_{m=0}^{\infty}E_{m}^{\left(s\right)}\frac{t^{m}}{m!}\right)\nonumber \\
 & =\sum_{n=0}^{\infty}\left(\sum_{l=0}^{n}\dbinom{n}{l}B_{l}^{\left(r\right)}\left(x\right)E_{n-l}\right)\frac{t^{n}}{n!},\nonumber
\end{align}
and
\begin{align}
\sum_{m=0}^{\infty}DC_{n}^{\left(r,s\right)}\left(x\right)\frac{\left(e^{t}-1\right)^{m}}{m!} & =\sum_{m=0}^{\infty}DC_{m}^{\left(r,s\right)}\left(x\right)\sum_{n=m}^{\infty}S_{2}\left(n,m\right)\frac{t^{n}}{n!}\label{eq:33}\\
 & =\sum_{n=0}^{\infty}\left(\sum_{m=0}^{n}DC_{m}^{\left(r,s\right)}\left(x\right)S_{2}\left(n,m\right)\right)\frac{t^{n}}{n!}.\nonumber
\end{align}

Therefore, by (\ref{eq:32}) and (\ref{eq:33}), we get
\begin{equation}
\sum_{m=0}^{n}DC_{m}^{\left(r,s\right)}\left(x\right)S_{2}\left(m,n\right)=\sum_{l=0}^{n}\dbinom{n}{l}B_{l}^{\left(r\right)}\left(x\right)E_{n-l},\label{eq:34}
\end{equation}
where $S_{2}\left(n,m\right)$ is the Stirling number of the second
kind.

Finally, we consider the Cauchy-Changhee mixed-type polynomials of order $(r,s)$ as
follows :
\begin{equation}
CC_{n}^{\left(r,s\right)}\left(x\right)=\int_{\Zp}\cdots\int_{\Zp}C_{n}^{\left(r\right)}\left(x+y_{1}+\cdots+y_{s}\right)d\mu_{-1}\left(y_{1}\right)\cdots d\mu_{-1}\left(y_{s}\right),\label{eq:35}
\end{equation}
where $n\ge0$.

By (\ref{eq:35}), we see that the generating function of $CC_{n}^{\left(r,s\right)}\left(x\right)$
are given by
\begin{align}
 & \sum_{n=0}^{\infty}CC_{n}^{\left(r,s\right)}\left(x\right)\frac{t^{n}}{n!}\label{eq:36}\\
= & \int_{\Zp}\cdots\int_{\Zp}\sum_{n=0}^{\infty}C_{n}^{\left(r\right)}\left(x+y_{1}+\cdots+y_{s}\right)\frac{t^{n}}{n!}d\mu_{-1}\left(y_{1}\right)\cdots d\mu_{-1}\left(y_{s}\right)\nonumber \\
= & \left(\frac{t}{\log\left(1+t\right)}\right)^{r}\int_{\Zp}\cdots\int_{\Zp}\left(1+t\right)^{x+y_{1}+\cdots+y_{s}}d\mu_{-1}\left(y_{1}\right)\cdots d\mu_{-1}\left(y_{s}\right)\nonumber \\
= & \left(\frac{t}{\log\left(1+t\right)}\right)^{r}\left(\frac{2}{t+2}\right)^{s}\left(1+t\right)^{x}\nonumber \\
= & \left(\sum_{m=0}^{\infty}C_{m}^{\left(r\right)}\left(x\right)\frac{t^{m}}{m!}\right)\left(\sum_{l=0}^{\infty}Ch_{l}^{\left(s\right)}\frac{t^{l}}{l!}\right)\nonumber \\
= & \sum_{n=0}^{\infty}\left\{ \sum_{m=0}^{n}\dbinom{n}{m}C_{m}^{\left(r\right)}\left(x\right)Ch_{n-m}^{\left(s\right)}\right\} \frac{t^{n}}{n!}.\nonumber
\end{align}

Thus, by (\ref{eq:36}), we get
\begin{equation}
CC_{n}^{\left(r,s\right)}\left(x\right)=\sum_{m=0}^{n}\dbinom{n}{m}C_{m}^{\left(r\right)}\left(x\right)Ch_{n-m}^{\left(s\right)}.\label{eq:37}
\end{equation}

By replacing $t$ by $e^{t}-1$, we get
\begin{align}
\sum_{n=0}^{\infty}CC_{n}^{\left(r,s\right)}\left(x\right)\frac{\left(e^{t}-1\right)^{n}}{n!} & =\left(\frac{e^{t}-1}{t}\right)^{r}\left(\frac{2}{e^{t}+1}\right)^{s}e^{xt}\label{eq:38}\\
 & =\left(\sum_{l=0}^{\infty}\frac{S_{2}\left(l+r,l\right)}{\binom{l+r}{r}}\frac{t^{l}}{l!}\right)\left(\sum_{m=0}^{\infty}E_{m}^{\left(s\right)}\left(x\right)\frac{t^{m}}{m!}\right)\nonumber \\
 & =\sum_{n=0}^{\infty}\left(\sum_{l=0}^{n}\frac{S_{2}\left(l+r,l\right)E_{n-l}^{\left(s\right)}\left(x\right)}{\binom{l+r}{l}}\dbinom{n}{l}\right)\frac{t^{n}}{n!},\nonumber
\end{align}
and
\begin{align}
\sum_{l=0}^{\infty}CC_{l}^{\left(r,s\right)}\left(x\right)\frac{\left(e^{t}-1\right)^{l}}{l!} & =\sum_{l=0}^{\infty}CC_{l}^{\left(r,s\right)}\left(x\right)\sum_{n=l}^{\infty}S_{2}\left(n,l\right)\frac{t^{n}}{n!}\label{eq:39}\\
 & =\sum_{n=0}^{\infty}\left(\sum_{l=0}^{n}CC_{l}^{\left(r,s\right)}\left(x\right)S_{2}\left(n,l\right)\right)\frac{t^{n}}{n!}.\nonumber
\end{align}

Therefore, by (\ref{eq:38}) and (\ref{eq:39}), we obtain the following
identities.

\[
\sum_{l=0}^{n}CC_{l}^{\left(r,s\right)}\left(x\right)S_{2}\left(n,l\right)=\sum_{l=0}^{n}\frac{\binom{n}{l}}{\binom{l+r}{l}}S_{2}\left(l+r,l\right)E_{n-l}^{\left(s\right)}\left(x\right),
\]
where $n\ge0$.

\providecommand{\bysame}{\leavevmode\hbox to3em{\hrulefill}\thinspace}
\providecommand{\MR}{\relax\ifhmode\unskip\space\fi MR }
\providecommand{\MRhref}[2]{%
  \href{http://www.ams.org/mathscinet-getitem?mr=#1}{#2}
}
\providecommand{\href}[2]{#2}

\end{document}